\documentclass[12pt]{amsart}
\usepackage{amsmath,amssymb,amsthm,amscd}

\usepackage{graphicx}

\newtheorem{formula}{}[section]
\newtheorem{proposition}[formula]{Proposition}
\newtheorem{corollary}[formula]{Corollary}
\newtheorem{lemma}[formula]{Lemma}
\newtheorem{theorem}[formula]{Theorem}

\theoremstyle{definition}
\newtheorem{definition}[formula]{Definition}
\newtheorem{construction}[formula]{Construction}

\theoremstyle{remark}
\newtheorem{remark}{Remark}
\newtheorem{problem}{Problem}

\begin{document}

\title[$B$-rigidity of ideal almost Pogorelov polytopes]{$B$-rigidity of ideal almost Pogorelov polytopes}
\author[N.Yu.~Erokhovets]{Nikolai~Erokhovets}
\address{Department of Mechanics and Mathematics, Moscow State University, 119992 Moscow, Russia}
\email{erochovetsn@hotmail.com}

\def\sgn{\mathrm{sgn}\,}
\def\bideg{\mathrm{bideg}\,}
\def\tdeg{\mathrm{tdeg}\,}
\def\sdeg{\mathrm{sdeg}\,}
\def\grad{\mathrm{grad}\,}
\def\ch{\mathrm{ch}\,}
\def\sh{\mathrm{sh}\,}
\def\th{\mathrm{th}\,}

\def\mod{\mathrm{mod}\,}
\def\In{\mathrm{In}\,}
\def\Im{\mathrm{Im}\,}
\def\Ker{\mathrm{Ker}\,}
\def\Hom{\mathrm{Hom}\,}
\def\Tor{\mathrm{Tor}\,}
\def\rk{\mathrm{rk}\,}
\def\codim{\mathrm{codim}\,}

\def\ko{{\mathbf k}}
\def\sk{\mathrm{sk}\,}
\def\RC{\mathrm{RC}\,}
\def\gr{\mathrm{gr}\,}

\def\R{{\mathbb R}}
\def\C{{\mathbb C}}
\def\Z{{\mathbb Z}}
\def\A{{\mathcal A}}
\def\B{{\mathcal B}}
\def\K{{\mathcal K}}
\def\M{{\mathcal M}}
\def\N{{\mathcal N}}
\def\E{{\mathcal E}}
\def\G{{\mathcal G}}
\def\D{{\mathcal D}}
\def\F{{\mathcal F}}
\def\L{{\mathcal L}}
\def\V{{\mathcal V}}
\def\H{{\mathcal H}}

\thanks{The research was supported by the RFBR grant No 20-01-00675}





\begin{abstract}
Toric topology assigns to each $n$-dimensional combinatorial simple convex polytope $P$ with $m$ facets an 
$(m+n)$-dimensional moment-angle manifold $\mathcal{Z}_P$ with an action of a compact torus $T^m$ such that 
$\mathcal{Z}_P/T^m$ is a convex polytope of combinatorial type $P$.  A simple  $n$-polytope is called $B$-rigid,
if any isomorphism of graded rings $H^*(\mathcal{Z}_P,\mathbb Z)= H^*(\mathcal{Z}_Q,\mathbb Z)$ 
for a simple $n$-polytope $Q$ implies that $P$ and $Q$ are combinatorially equivalent. 
An ideal almost Pogorelov polytope is a combinatorial $3$-polytope 
obtained by cutting off all the ideal vertices of an ideal right-angled polytope in the Lobachevsky (hyperbolic) space 
$\mathbb L^3$. These polytopes are exactly the polytopes obtained from any, not necessarily simple, convex $3$-polytopes
by cutting off all the vertices followed by cutting off all the "old" edges. The boundary of the dual polytope
is the barycentric subdivision of the boundary of the old polytope (and also of its dual polytope). We prove 
that any ideal almost Pogorelov polytope is $B$-rigid. 
This produces three cohomologically rigid families of manifolds over ideal almost Pogorelov manifolds:
moment-angle manifolds, canonical $6$-dimensional quasitoric manifolds and canonical $3$-dimensional small covers,
which are "pullbacks from the linear model". 
\end{abstract}
\maketitle

\setcounter{section}{0}

\section*{Introduction}
For an introduction to the polytope theory we recommend \cite{Z14}. For brevity, by an {\it $n$-polytope} we
call a convex $n$-dimensional polytope.  All the details 
concerning basic notions from toric topology can be found in \cite{BP15} and \cite{BE17S}. 
Here we give a brief description of these notions. Let $S^1=\{z\in\mathbb C\colon |z|=1\}$, and $T^m=(S^1)^m$ be a 
compact torus.
Toric topology assigns to each simple $n$-polytope $P$  with $m$ faces $F_1,\dots,F_m$ an $(m+n)$-dimensional 
{\it moment-angle manifold} 
$$
\mathcal{Z}_P= T^m\times P^n/\sim,\text{ where }(t_1,p_1)\sim(t_2,p_2)\text{ if and only if }p_1=p_2,\text{ and }t_1t_2^{-1}\in T^{G(p_1)},
$$ 
and $T^{G(p)}=\{(t_1,\dots,t_m)\in T^m\colon t_i=1\text{ for }F_i\not\ni p\}$.  
There is an action of $T^m$ on $\mathcal{Z}_P$ induced from its action on the first factor, and $\mathcal{Z}_P/T^m=P$. 
It can be shown that topological type of  $\mathcal{Z}_P$ depends only on combinatorial type of $P$ and $\mathcal{Z}_P$ has a smooth structure such that the action of $T^m$ is smooth.

The mapping $\Lambda\colon \{F_1,\dots,F_m\}\to \mathbb Z^n$ such that for any vertex $v=F_{i_1}\cap\dots\cap F_{i_n}$
the vectors $\Lambda(F_{i_1})$, $\dots$, $\Lambda(F_{i_n})$ form a basis in $\mathbb Z^n$
is called {\it characteristic}. It induces a mapping
$$
\varphi_{\Lambda}\colon T^m\to T^n,\quad \varphi_{\Lambda}(t_1,\dots,t_m)=(t_1^{\Lambda_{1,1}}t_2^{\Lambda_{1,2}}\dots t_m^{\Lambda_{1,m}},\dots,t_1^{\Lambda_{n,1}}
t_2^{\Lambda_{n,2}}\dots t_m^{\Lambda_{n,m}}),
$$
where $S^1=\{z\in\mathbb C\colon |z|=1\}$, and $\Lambda(F_i)=(\Lambda_{1,i},\dots,\Lambda_{n,i})\in\mathbb Z^n$.
It can be shown that the subgroup $\Ker \varphi_{\Lambda}$ acts on $\mathcal {Z}_P$ freely. The quotient
space is known as a {\it quasitoric manifold}: 
$$
M(P,\Lambda)=T^n\times P^n/\sim,
$$
where $(t_1,p_1)\sim(t_2,p_2)$ if and only if $p_1=p_2$,  and $t_1t_2^{-1}\in \varphi_{\Lambda}(T^{G(p_1)})$.
The quasitoric manifold has an action of $T^n=T^m/\Ker\varphi_{\Lambda}$, and $M(P,\Lambda)/T^n=P$.
These manifolds first appeared in \cite{DJ91} as topological analogs of toric manifolds from algebraic geometry.

Similarly, for a {\it $\mathbb Z_2$-characteristic mapping} $\Lambda_2\colon \{F_1,\dots,F_m\}\to\mathbb Z_2^n=
(\mathbb Z/2\mathbb Z)^n$ such that for any vertex $v=F_{i_1}\cap\dots\cap F_{i_n}$
the vectors $\Lambda_2(F_{i_1})$, $\dots$, $\Lambda_2(F_{i_n})$ form a basis in $\mathbb Z_2^n$ there is a {\it small cover}
$$
R(P,\Lambda_2)=\mathbb Z_2^n\times P^n/\sim,
$$
where $(s_1,p_1)\sim(s_2,p_2)$ if and only if $p_1=p_2$,  and $s_1-s_2\in \langle \Lambda_2(F_i)\rangle_{F_i\ni p_1}$.
It can be shown that any small cover is a smooth $n$-dimensional manifold glued from $2^n$ copies of $P$. Moreover, 
if $P$ and $Q$ are combinatorially equivalent, then the corresponding small covers are diffeomorphic. The small cover
has a canonical action of $\mathbb Z_2^n$ such that $R(P,\Lambda_2)/\mathbb Z_2^n=P$.
\begin{definition}
A simple $n$-polytope $P$ is called {\it $B$-rigid}, if for any simple $n$-polytope $Q$ any isomorphism of graded rings 
$H^*(\mathcal{Z}_P,\mathbb Z)\simeq H^*(\mathcal{Z}_Q,\mathbb Z)$ implies that $P$ and $Q$ are combinatorially equivalent.

A simple $n$-polytope $P$ is called {\it $C$-rigid}, if it 
admits a characteristic mapping, and for any simple $n$-polytope $Q$ any isomorphism of 
graded rings $H^*(M(P,\Lambda),\mathbb Z)\simeq H^*(M(Q,\Lambda'),\mathbb Z)$ implies that $P$ and $Q$ 
are combinatorially equivalent.
\end{definition}
It is known that $B$-rigid polytope is $C$-rigid \cite{CPS10} (see also \cite{BEMPP17}). The converse is not
true \cite{CP19}. 

In what follows, unless otherwise specified, by a {\it polytope} we mean a  class
of combinatorially equivalent $3$-polytopes. By {\it faces} we call its facets.

For any simple $3$-polytope $P$ by the Four colour theorem 
there is a colouring $c\colon \{F_1,\dots,F_m\}\to\{1,2,3,4\}$ such that adjacent faces have different colours. 
Each colouring $c$ induces a characteristic function $\Lambda_c$ by the following rule.
Let $\boldsymbol{e}_1$, $\boldsymbol{e}_2$, and $\boldsymbol{e}_3$ be a basis in $\mathbb Z^3$, and 
$\boldsymbol{e}_4=\boldsymbol{e}_1+\boldsymbol{e}_2+\boldsymbol{e}_3$. Set $\Lambda_c(F_i)=\boldsymbol{e}_{c(i)}$.
Then  $\Lambda_c$ is a characteristic function, since any three vectors from the set 
$\{\boldsymbol{e}_1, \boldsymbol{e}_2, \boldsymbol{e}_3, \boldsymbol{e}_4\}$ form a basis in $\mathbb Z^3$.

\begin{definition}
A {\it $k$-belt} is a cyclic sequence of $k$ faces with the property that faces are adjacent if and only if they follow each other, and no three faces have a common vertex.  A $k$-belt is {\it trivial}, if it surrounds a face. 
\end{definition}

\begin{definition}
A simple $n$-polytope $P$ is called {\it flag}, if any set of its pairwise intersecting faces $F_{i_1},\dots,F_{i_k}$ has a 
non-empty intersection $F_{i_1}\cap\dots\cap F_{i_k}\ne \varnothing$. 
\end{definition}

It can be shown (see \cite{BE15}) that a simple $3$-polytope 
$P$ is flag if and only if $P\ne\Delta^3$, and $P$ has no $3$-belts.

Results by  A.V.~Pogorelov  \cite{P67} and E.M.~Andreev \cite{A70a} imply that a simple $3$-polytope $P$
can be realized in the Lobachevsky (hyperbolic) space $\mathbb L^3$ as a bounded polytope with right dihedral angles
if and only if $P$ is different from the simplex $\Delta^3$ and has no $3$- and $4$-belts. Moreover, a realization is unique up to 
isometries. Such polytopes are called
{\it Pogorelov polytopes}. We denote their family $\mathcal{P}_{Pog}$. 
It follows from \cite{D98,D03} (see also \cite{BE17I}) that  $\mathcal{P}_{Pog}$ contains  {\it fullerenes}, that is simple $3$-polytopes with only $5$- and $6$-gonal. Mathematical fullerenes model spherical carbon atoms. In 1996   R.~Curl, H.~Kroto, and R.~Smalley obtained the Nobel Prize in chemistry ``for their discovery of fullerenes''. They  synthesized  {\it Buckminsterfullerene} $C_{60}$, which has the form of the truncated icosahedron (and also of a soccer ball). W.P.~Thurston \cite{T98} built a parametrisation of the fullerene family, which implies that the number of fullerenes with  $n$ carbon atoms grows like~$n^9$ when $n$ tends to infinity. 

In \cite{FMW15} F.~Fan, J.~Ma, and X.~Wang proved that any Pogorelov polytope is $B$-rigid. In particular, it is $C$-rigid.

\begin{remark}
The notions of $B$- and $C$-rigidity can be defined for any field $\mathbb F$ instead of $\mathbb Z$. Then
a polytope $B$-rigid over $\mathbb F$ is $C$-rigid over $\mathbb F$. In fact, in paper \cite{FMW15} it is proved that
any Pogorelov polytope is $B$-rigid over $\mathbb Z$ or any field $\mathbb F$. 
\end{remark}
\begin{remark}
In \cite{B17} F.~Bosio presented a construction of flag $3$-polytopes, which are not $B$-rigid.
\end{remark}

In \cite{BEMPP17} for Pogorelov polytopes a fact stronger than $C$-rigidity was proved. Namely, 
it is known (see \cite{DJ91,BP15,BEMPP17}) that the manifolds $M(P,\Lambda)$ and $M(Q,\Lambda')$ over two simple 
$n$-polytopes $P$ and $Q$ are weakly equivariantly diffeomorphic, that 
is there is a diffeomorphism $f\colon M(P,\Lambda)\to M(Q,\Lambda')$ and an automorphism $\psi\colon T^n\to T^n$ such that 
$f(t\cdot x)=\psi(t)\cdot f(x)$ for all $t\in T^n$ and $x\in M(P,\Lambda)$, if and only if there is a combinatorial
equivalence $\varphi\colon P\to Q$ and a change of basis $C\in GL(n, \mathbb Z)$ such that 
$\Lambda(\varphi(F_i))=\pm C\cdot \Lambda(F_i)$ for $i=1,\dots, m$. Such pairs $(P,\Lambda)$ and $(Q,\Lambda')$ are 
called {\it equivalent}. 

\begin{theorem}[\cite{BEMPP17}]\label{C6th} Let $P,Q\in \mathcal{P}_{Pog}$. Then there is an isomorphism of graded rings 
$H^*(M(P,\Lambda),\mathbb Z)\simeq H^*(M(Q,\Lambda'),\mathbb Z)$ if and only 
if the pairs $(P,\Lambda)$ and $(Q, \Lambda')$ are equivalent. 
\end{theorem}
Moreover, in \cite{BP16} it was proved that the pairs 
$(P,\Lambda_c)$ and $(Q, \Lambda'_{c'})$ induced by colourings are equivalent if and only if there is a combinatorial 
equivalence $\varphi\colon P\to Q$ and a permutation $\sigma\in S_4$ such that $c'(\varphi(F_i))=\sigma\circ c(F_i)$ for all $i$.

A similar result was proved for small covers. Namely, any characteristic function $\Lambda$ induces a 
$\mathbb Z_2$-characteristic function $\Lambda_2$ simply by taking all the coefficients modulo $2$. In dimension $3$
the converse is also true. The inclusion of sets $\mathbb Z_2^3\subset \mathbb Z^3$ lifts a $\mathbb Z_2$-characteristic
function $\Lambda_2$ to a characteristic function $\Lambda$, since for any matrix of size $3\times 3$ 
with entries $0$ and $1$, which is non-degenerate over $\mathbb Z_2$, the determinant is $\pm1$ in $\mathbb Z$.
This is not valid for dimensions greater than $3$. It is known \cite{DJ91} that 
$$
H^*(R(P,\Lambda_2),\mathbb Z_2)=\mathbb Z_2[v_1,\dots,v_m]/(I_P+J_{\Lambda_2}),
$$
where each variable $v_i$ has degree $1$, the ideal $I_P$ is generated by monomials $v_{i_1}\dots v_{i_k}$ such that $F_{i_1}\cap\dots\cap F_{i_k}=\varnothing$,
and the ideal $J_{\Lambda_2}$ is generated by linear forms $\Lambda_{1,1}v_1+\dots+\Lambda_{1,m}v_m$, $\dots$,
$\Lambda_{n,1}v_1+\dots+\Lambda_{n,m}v_m$. Similarly for a quasitoric manifold $M(P,\Lambda)$, and the coefficient ring $R$, which is $\mathbb Z$ or $\mathbb  Z_2$
we have \cite{DJ91}:
$$
H^*(M(P,\Lambda), R)= R[v_1,\dots,v_m]/(I_P+J_{\Lambda_2}),
$$
where each variable $v_i$ has degree $2$.

Thus, any isomorphism of graded rings $\varphi\colon H^*(R(P,\Lambda_2),\mathbb Z_2)\to H^*(R(Q,\Lambda_2'),\mathbb Z_2)$
induces an isomorphism for $M(P,\Lambda)$ and $M(Q,\Lambda')$. In particular, if $P$ and $Q$ are Pogorelov
polytopes, then they are combinatorially equivalent. Moreover, it was proved in \cite{BEMPP17}
that the pairs $(P,\Lambda_2)$ and $(Q,\Lambda_2')$ in this case are {\it $\mathbb Z_2$-equivalent}, that is there is a
a combinatorial equivalence $\varphi\colon P\to Q$ and a change of basis $C\in GL(3, \mathbb Z_2)$ such that 
$\Lambda(\varphi(F_i))=C\cdot \Lambda(F_i)$ for $i=1,\dots, m$.

Recent results of this type in the context of toric varieties for $n$-dimensional cubes see in \cite{CLMP20}.

\begin{construction}\label{Vconstr}
In \cite{V87} A.Yu.~Vesnin introduced a construction of a  $3$-dimensional compact hyperbolic manifold 
$N(P,\Lambda_2)$ corresponding to a Pogorelov polytope $P$ and a $\mathbb Z_2$-characteristic mapping $\Lambda_2$. 
Namely, for a bounded right-angled polytope $P\subset\mathbb L^3$ 
there is a right-angled Coxeter  group $G(P)$ generated by reflections in faces of $P$. It it known that 
$$
G(P)=\langle \rho_1,\dots,\rho_m\rangle/(\rho_1^2,\dots,\rho_m^2, \rho_i\rho_j=\rho_j\rho_i
\text{ for all }F_i\cap F_j\ne\varnothing),
$$
where $\rho_i$ is the reflection in the face $F_i$. The group $G(P)$ acts on $\mathbb L^3$ discretely (see details in \cite{VS88}),
and $P$ is a fundamental domain of $G(P)$, that is the interiors of the polytopes $\{gP\}_{g\in G(P)}$ do not intersect. 
Moreover, the stabiliser of a point $x\in P$ is generated by reflections in faces containing it.   

The mapping $\Lambda_2$ defines a homomorphism $\varphi\colon G(P)\to \mathbb Z_2^3$
by the rule $\varphi(\rho_i)=\Lambda_2(F_i)$. It can be shown that $\Ker\varphi$ acts on $\mathbb L^3$ freely. 
Then $N(P,\Lambda_2)=\mathbb L^3/\Ker\varphi$ is a compact hyperbolic
manifold glued of $8$ copies of $P$ along faces. It is easy to see that $N(P,\Lambda_2)$ is homeomorphic to a small
cover $R(P,\Lambda_2)$, and the homeomorphism is given by the mapping $(t,p)\to \varphi^{-1}(t)\cdot p$. 
\end{construction}

In this paper we consider another family $\mathcal{P}_{aPog}$ of $3$-polytopes. It consists of simple $3$-polytopes 
$P\ne \Delta^3$ without $3$-belts such that any $4$-belt surrounds a face. We call them {\it almost Pogorelov polytopes}. 
The combinatorics and hyperbolic geometry of the family $\mathcal{P}_{aPog}$ was studied in \cite{E19}. 

In 2019 T.E.~Panov remarked that results by E.~M.~Andreev \cite{A70a, A70b} should imply that almost Pogorelov 
polytopes correspond to right-angled polytopes of finite volume in  $\mathbb L^3$. 
Such polytopes may have $4$-valent vertices on the absolute, while all proper vertices have valency $3$. 
It was proved in \cite{E19} that  cutting of $4$-valent vertices defines a bijection between classes of congruence of right-angled
polytopes of finite volume in $\mathbb {L}^3$ and almost Pogorelov polytopes different from the cube $I^3$ 
and the pentagonal prism $M_5\times I$. Moreover, it induces the bijection between the 
ideal vertices of the right-angled polytope and the 
quadrangles of the corresponding almost Pogorelov polytope. The polytopes
$I^3$ and  $M_5\times I$ are the only almost Pogorelov polytopes with adjacent quadrangles. The cube has $m=6$ faces. The pentagonal prism has
$7$ faces. For $m=8$ there are no almost Pogorelov polytopes. For $m=9$ there is a unique almost Pogorelov polytope --
the $3$-dimensional associahedron $As^3$. It can be realized as a cube with three non-adjacent pairwise orthogonal edges cut.
It follows from \cite{B74} that all the other polytopes in $\mathcal{P}_{aPog}$ can be obtained from $As^3$ by a sequence of
operations of cutting off an edge not lying in quadrangles, and cutting off a pair of adjacent edges of a face with 
at least $6$ sides by one plane. Moreover, these operations preserve the family of almost Pogorelov polytopes.

\begin{figure}
  \centering
  \includegraphics[height=4cm]{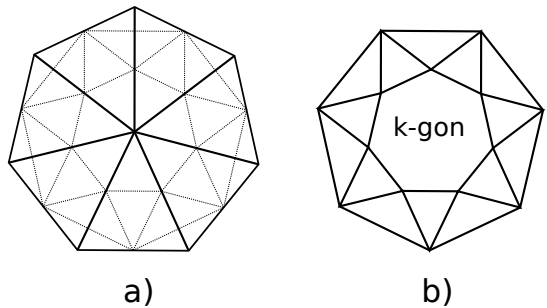}
  \caption{a) Medial graph of the $k$-gonal pyramid; b) The $k$-antiprism}
  \label{k-antiprism}
\end{figure}

\begin{definition}
A polytope in $\mathbb{L}^3$ is called {\it ideal}, if all its vertices lie on the absolute (are {\it ideal points}). 
An ideal polytope has a finite volume. 
Its congruence class is uniquely defined by the combinatorial type \cite{R94} (see also \cite{R96}).
We will call almost Pogorelov polytopes corresponding to ideal right-angled polytopes {\it ideal almost Pogorelov polytopes} 
and denote their family $\mathcal{P}_{IPog}$.
\end{definition}
It is a classical fact (see, for example \cite[Section 13.6]{T02} and \cite[Theorem 5 (iv)]{BGGMTW05}) 
that graphs of ideal right-angled $3$-polytopes  are exactly {\it medial graphs} 
of $3$-polytopes. Vertices of a medial graph correspond to edges of a polytope, and edges -- to pairs of edges adjacent 
in a face. This correspondence plays a fundamental role in the well-known Koebe-Andreev-Thurston theorem 
(see \cite{Z14,T02,BS04,S92}): {\it any $3$-dimensional combinatorial polytope $P$  has a geometric realization 
in the Euclidean space $\mathbb R^3$ such that all its edges are tangent to a given sphere}.  
A medial graph of the $k$-pyramid is known as $k$-antiprism, see Fig. \ref{k-antiprism} b). It follows from the definition
that ideal almost Pogorelov polytope $P$ corresponding to a polytope $Q$ can be obtained from it by cutting off all the vertices
followed by cutting off all the "old" edges. Such polytopes appeared in  \cite[Example 1.15(1)]{DJ91} 
as polytopes giving examples for quasitoric manifolds and small covers, which are "pullbacks from the linear model". 
We prove that for $3$-polytopes the families of these manifolds are cohomologically rigid 
(see Theorem \ref{FCth}). The triangulation  $\partial P^*$ is the barycentric subdivision of $\partial Q^*$. 
Such triangulations appeared also in \cite[Example 6.3(2)]{FMW20} as examples of flag triangulations such that any $4$-circuit 
full subcomplex is simple. Also they appeared in \cite{CK11} in study of {\it combinatorial rigidity} 
(a simple $n$-polytope $P$ is called {\it combinatorially rigid}, if there are no other simple $n$-polytopes with the same bigraded Betti numbers $\beta^{-i,2j}(\mathcal{Z}_P)$).

\begin{figure}
  \centering
  \includegraphics[width=200pt]{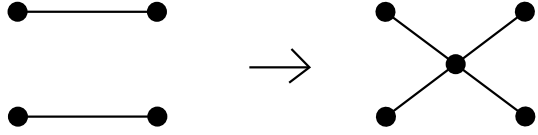}
  \caption{An edge-twist}
  \label{E-twist}
\end{figure}

\begin{figure}
\begin{center}
\includegraphics[height=3cm]{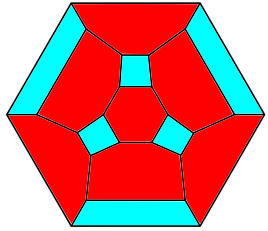}
\caption{Three-dimensional permutohedron $Pe^3$}\label{Pe3}
\end{center}
\end{figure}

An operation of an {\it edge-twist} is drawn on Fig. \ref{E-twist}. Two edges on the left lie in the same face and are disjoint. Let us call an edge-twist {\it restricted}, if both edges are adjacent to an edge of the same face. It follows from \cite{V17,BGGMTW05, E19} that  a polytope is realizable as an ideal right-angled polytope if and only if either it is a $k$-antiprism, $k\geqslant 3$, or it can be obtained from the $4$-antiprism by a sequence of restricted edge-twists.  
Thus, the simplest ideal right-angled polytope is a $3$-antiprism, which coincides with the octahedron.
The corresponding ideal almost Pogorelov polytope is known as $3$-dimensional permutohedron $Pe^3$, see Fig. \ref{Pe3}. This is a unique ideal almost Pogorelov polytope with minimal number $m$ of faces ($14$ faces). It is easy to see that 
for ideal almost Pogorelov polytopes $m=2(p_4+1)$, where $p_4$ is the number of quadrangles. 
Recent results on volumes of ideal right-angled polytopes see in \cite{VE20}.

In the paper \cite{E20} the technique from \cite{FMW15} (see also \cite{BE17S})  was generalized to the family 
$\mathcal{P}_{aPog}$. In particular, it was proved that if there is an isomorphism of graded rings 
$H^*(\mathcal{Z}_P,\mathbb Z)\simeq H^*(\mathcal{Z}_Q,\mathbb Z)$, where $P$ is an almost Pogorelov or 
an ideal almost Pogorelov polytope,
and $Q$ is any simple $3$-polytope, then 
$Q$ is an almost Pogorelov \cite[Theorem 3.6]{E20} or an ideal almost Pogorelov polytope 
\cite[Corollary 3.8]{E20} respectively.

Main result of our paper is 
\begin{theorem}\label{IBthe}
Any ideal almost Pogorelov polytope is $B$-rigid.
\end{theorem}
The proof follows from results in \cite{E20}. We compare it with the proof of $B$-rigidity of Pogorelov polytopes 
from \cite{FMW15}. 

In Theorem \ref{FCth} we show how Theorem \ref{IBthe} gives 
three cohomologically rigid families of manifolds over ideal almost Pogorelov manifolds:
moment-angle manifolds, canonical $6$-dimensional quasitoric manifolds and canonical $3$-dimensional small covers,
which are "pullbacks from the linear model" from \cite[Example 1.15(1)]{DJ91}.  

Below is the plan of the paper.

In Section \ref{Cohsec} we give a brief description of the cohomology ring
of a moment-angle manifold.

In Section \ref{SCCsec} we describe an important tool in the study of families if $3$-polytopes:
a separable circuit condition.

In Section \ref{Brsec} we describe a notion of a $B$-rigid property, subset, a collection of subsets, or a number.
We explain, why the properties to be a flag, Pogorelov and almost Pogorelov polytopes are
$B$-rigid. Also there are several interesting open problems in this section. 

In Section \ref{3dsec} we describe $B$-rigid subsets in $H^3(\mathcal{Z}_P)$ for Pogorelov 
and almost Pogorelov polytopes.

In Section \ref{Bcsec} we explain, which sets of 
elements in $H^{k+2}(\mathcal{Z}_P)$ corresponding to $k$-belts are $B$-rigid
for Pogorelov and almost Pogorelov polytopes

In Section \ref{Btsec} we deal with belts surrounding faces. It turns out that for Pogorelov polytopes or ideal almost
Pogorelov polytopes $P,Q$ the elements in $H^*(\mathcal{Z}_P)$ corresponding to belts around faces under each 
isomorphism of graded rings $H^*(\mathcal{Z}_P)\to H^*(\mathcal{Z}_Q)$ are mapped to the analogous elements for $Q$.
This induces a bijection between the sets of faces of $P$ and $Q$. 

In Section \ref{Basec} we discuss that under this bijection adjacent faces are mapped to adjacent faces.
In particular, there is a combinatorial equivalence between $P$ and $Q$. This proves Theorem \ref{IBthe}.

In Section \ref{Newsec} we show  (Theorem \ref{FCth}) how Theorem \ref{IBthe} produces
three cohomologically rigid families of manifolds over ideal almost Pogorelov manifolds:
moment-angle manifolds, canonical $6$-dimensional quasitoric manifolds and canonical $3$-dimensional small covers,
which are "pullbacks from the linear model" from \cite[Example 1.15(1)]{DJ91}.  

In Section \ref{Hsec} we study hyperbolic geometry of the $3$-dimensional manifolds corresponding to 
ideal almost Pogorelov polytopes.

\begin{remark}
In \cite[Theorems 3.12 and 3.14]{FMW20} F.~Fan, J.~Ma, and X.~Wang generalized 
results from \cite{FMW15} and \cite{BEMPP17} to flag polytopes of dimension higher than $3$ 
satisfying a generalization of the separable circuit condition  (see Section \ref{SCCsec}).  
It is announced without any details, that in the next paper \cite{FMW20b} the authors will prove 
$B$-rigidity of any almost Pogorelov polytope and an analog of Theorem \ref{C6th} for them. 
But their \cite[Definition 2.16]{FMW20} of $B$-rigidity assumes isomorphism of bigraded rings, which is a more restrictive condition.
\end{remark}

\section{Cohomology ring of a moment-angle manifold of a simple $3$-polytope}\label{Cohsec}
Details on cohomology of a moment-angle manifolds see in \cite{BP15, BE17S}.
If not specified, we study cohomology over $\mathbb Z$ and omit the coefficient ring.

The ring $H^*(\mathcal{Z}_P)$  has a multigraded structure:
$$
H^*(\mathcal{Z}_P)=\bigoplus\limits_{i\geqslant0,\omega\subset[m]}H^{-i,2\omega}(\mathcal{Z}_P),
\text{ where }H^{-i,2\omega} (\mathcal{Z}_P)\subset H^{2|\omega|-i} (\mathcal{Z}_P),
$$
and $[m]=\{1,\dots,m\}$. There is a canonical isomorphism 
$H^{-i,2\omega}(\mathcal{Z}_P)\simeq \widetilde{H}^{|\omega|-i-1}(P_{\omega})$,
where $P_{\omega}=\bigcup_{i\in \omega}F_i$, and $\widetilde{H}^{-1}(\varnothing)=\mathbb Z$.

The multiplication between components is nonzero, only if $\omega_1\cap \omega_2=\varnothing$. The mapping
$$
H^{-i,2\omega_1}(\mathcal{Z}_P)\otimes H^{-j,2\omega_2}(\mathcal{Z}_P)\to 
H^{-(i+j),2(\omega_1\sqcup\omega_2)}(\mathcal{Z}_P)
$$
via the Poincare-Lefschets duality $H^i(P_{\omega})\simeq H_{n-1-i}(P_{\omega},\partial P_{\omega})$ 
up to signs is induced by intersection of faces of $P$.

For any simple $3$-polytope $P$ and $\omega\ne\varnothing$ the set 
$P_{\omega}$ is a $2$-dimensional manifold, perhaps with a boundary. Therefore, $\widetilde{H}^k(P_{\omega})$
is nonzero only for $(k,\omega)\in \{(-1,\varnothing),(0,*),(1,*),(2,[m])\}$. In particular, $P$ has no torsion, and
there is an isomorphism $H^*(\mathcal{Z}_P,\mathbb Q)\simeq H^*(\mathcal{Z}_P)\otimes\mathbb Q$ and 
the embedding $H^*(\mathcal{Z}_P)\subset H^*(\mathcal{Z}_P)\otimes\mathbb Q$. For polytopes $P$ and $Q$ the 
isomorphism $H^*(\mathcal{Z}_P)\simeq H^*(\mathcal{Z}_Q)$ implies the isomorphism over $\mathbb Q$. 
For cohomology over $\mathbb Q$ or any field $\mathbb F$ 
all  theorems about the structure of $H^*(\mathcal{Z}_P,\mathbb Q)$ are still valid.

 There is a multigraded Poincare duality, which means that the bilinear form  
$$
H^{-i,2\omega}(\mathcal{Z}_P)\times H^{-(m-3-i),2([m]\setminus\omega)}(\mathcal{Z}_P)\to H^{-(m-3),2[m]}(\mathcal{Z}_P)=\mathbb Z
$$
has in some basis matrix with determinant $\pm1$.  We have 
\begin{gather*}
H^0(\mathcal{Z}_P)=\widetilde{H}^{-1}(\varnothing)=\mathbb Z=\widetilde{H}^2(\partial P)=H^{m+3}(\mathcal{Z}_P);\\
H^1(\mathcal{Z}_P)=H^2(\mathcal{Z}_P)=0=H^{m+1}(\mathcal{Z}_P)=H^{m+2}(\mathcal{Z}_P);\\
H^k(\mathcal{Z}_P)=\bigoplus\limits_{|\omega|=k-1}\widetilde{H}^0(P_{\omega})\oplus\bigoplus\limits_{|\omega|=k-2}\widetilde{H}^1(P_{\omega}),\quad 3\leqslant k\leqslant m.
\end{gather*}

Nontrivial multiplication occurs only for
\begin{enumerate}
\item $\widetilde{H}^{-1}(P_{\varnothing})\otimes \widetilde{H}^{k}(P_{\omega})\to \widetilde{H}^{k}(P_{\omega})$. This 
corresponds to multiplication by $1$ in $H^*(\mathcal{Z}_P)$.
\item  $\widetilde{H}^{0}(P_{\omega})\otimes \widetilde{H}^{1}(P_{[m]\setminus\omega})\to \widetilde{H}^{2}(\partial P)$.
This corresponds to the Poincare duality pairing.
\item $\widetilde{H}^{0}(P_{\omega_1})\otimes \widetilde{H}^{0}(P_{\omega_2})\to \widetilde{H}^{1}(P_{\omega_1\sqcup\omega_2})$. 
\end{enumerate}

A crucial role in $B$-rigidity questions for $3$-polytopes is played by the following two types of classes 
in $H^*(\mathcal{Z}_P)$. Denote 
$$
N_2(P)=\{\{i,j\}\subset [m]\colon F_i\cap F_j=\varnothing\}.
$$
For any $\omega\in N_2(P)$ we have $\widetilde H^0(P_{\omega})=\mathbb Z$. Choose a generator 
in this group and denote it by $\widetilde{\omega}$. 
Then for any $3$-polytope $H^3(\mathcal{Z}_P)$ is a free abelian
group with the basis $\{\widetilde\omega\colon\omega\in N_2(P)\}$.

For any $k$-belt $\mathcal{B}_k$ we have $\widetilde{H}^1(P_{\omega})=\mathbb Z$ for
$\omega=\{i\colon F_i\in\mathcal{B}_k\}$. Choose a generator 
in this group and denote it by $\widetilde{\mathcal{B}_k}$.
Denote by $\boldsymbol{B}_k$ the subgroup in $H^{k+2}(\mathcal{Z}_P)$, $3\leqslant k\leqslant m-2$, with the basis 
$\{\widetilde{\mathcal{B}_k}\colon \mathcal{B}_k  \text{ is a $k$-belt}\}$.

\section{Separable circuit condition for families of $3$-polytopes}\label{SCCsec}
A very important tools for study of families of $3$-dimensional polytopes are given 
by so-called separable circuit conditions (SCC for short) \cite{FMW20}.
\begin{proposition}[SCC for flag polytopes \cite{FW15}]\label{SCCflag} A simple $3$-polytope
is flag if and only if  for  any three pairwise different faces $\{F_i,F_j,F_k\}$ with $F_i\cap F_j=\varnothing$ there exists 
an $l$-belt $\mathcal{B}_l$ such that $F_i,F_j\in \mathcal{B}_l$ and $F_k\notin \mathcal{B}_l$. 
\end{proposition}
Denote $|\mathcal{B}_l|=\bigcup_{F_i\in \mathcal{B}_l}F_i$.
\begin{proposition}[SCC for Pogorelov polytopes \cite{FMW15}]\label{SCCpog} A simple $3$-polytope
is Pogorelov if and only if  for any three pairwise different faces $\{F_i,F_j,F_k\}$ 
with $F_i\cap F_j=\varnothing$  there exists  an $l$-belt $\mathcal{B}_l$ such that 
$F_i,F_j\in \mathcal{B}_l$, $F_k\notin \mathcal{B}_l$, 
and $F_k$ does not intersect at least one of the two connected components of $|\mathcal{B}_l|\setminus(F_i\cup F_j)$.
\end{proposition}
\begin{remark}
In \cite{FMW20} SCC for Pogorelov polytopes was generalized to higher dimensions (we call it SCC' for short). In particular, 
the product of two flag polytopes with SCC' is also a flag polytope with SCC'. 
For a flag simple $n$-polytope with SCC' there is a construction \cite[Construction E.1]{FMW20} 
which under some assumptions gives a flag  $(n+1)$-polytope with SCC'.  
\end{remark}
Let $P_8$ be the cube with two non-adjacent orthogonal edges cut.
\begin{proposition}[SCC for almost Pogorelov polytopes \cite{FMW20,E20}]\label{SCCapog}
A simple $3$-polytope $P$ is an almost Pogorelov polytope  or the polytope $P_8$ if and only if for any three 
pairwise different faces $\{F_i,F_j,F_k\}$ such that $F_i\cap F_j=\varnothing$ an $l$-belt 
$\mathcal{B}_l$ such that $F_i,F_j\in \mathcal{B}_l$, $F_k\notin \mathcal{B}_l$, and $F_k$ does not intersect at least one of 
the two connected components of $|\mathcal{B}_l|\setminus(F_i\cup F_j)$ exists if and only if $F_k$ does not intersect 
quadrangles among the faces $F_i$ and $F_j$.
\end{proposition}
\begin{remark}
In \cite[Proposition F.1]{FMW20} this result was first proved in the dual setting for simplicial polytopes  in "only if" 
direction. 
In \cite{E20} this result was proved independently for simple polytopes in both directions.
\end{remark}

\section{$B$-rigid properties}\label{Brsec}
\begin{definition}
A property of an $n$-polytope $P$ is called {\it $B$-rigid}, if any isomorphism of graded rings $H^*(\mathcal{Z}_P,\mathbb Z)= H^*(\mathcal{Z}_Q,\mathbb Z)$ for a simple $n$-polytope $Q$ implies that it also has this property. 
\end{definition}

\begin{definition}\label{Bdef}
Let $\mathfrak{P}$ be some set of $3$-polytopes.

We call a number $n(P)$, a set $\mathfrak{S}_P\subset H^*(\mathcal{Z}_P)$ or a collection of such sets defined for any polytope $P\in\mathfrak{P}$ {\em $B$-rigid in the class} $\mathfrak{P}$ if for any isomorphism $\varphi$ of graded rings $H^*(\mathcal{Z}_P)\simeq H^*(\mathcal{Z}_Q)$, $P,Q\in \mathfrak{P}$, we have $n(P)=n(Q)$,  $\varphi(\mathfrak{S}_P)=\mathfrak{S}_Q$, or each set from the collection for $P$ is mapped bijectively to some set from the collection for $Q$ respectively.

To be short, $B$-rigidity in the class of all simple $3$-polytopes we call {\it $B$-rigidity}.
\end{definition}
\begin{proposition}
A property to be a flag $3$-polytope is $B$-rigid.
\end{proposition}
There are at least three different proofs of this fact. First, in \cite{FW15} it was proved that it is
equivalent to the fact that the ring $\widetilde{H}^*(\mathcal{Z}_P)/([\mathcal{Z}_P])$ is a (nonzero) indecomposable ring. Also it was proved that for flag $3$-polytopes 
$$
\widetilde{H}^1(P_{\omega})=\bigoplus_{\omega_1\sqcup\omega_2}\widetilde{H}^0(P_{\omega_1})\cdot\widetilde{H}^0(P_{\omega_2}).
$$
This result was based on SCC for flag polytopes (Proposition \ref{SCCflag}). 
Using the latter fact in \cite{BE17S} and \cite{BEMPP17} it was proved that a simple 
$3$-polytope $P\ne\Delta^3$ is flag if and only if 
$$
H^{m-2}(\mathcal{Z}_P)\subset(\widetilde{H}^*(\mathcal{Z}_P))^2.
$$
In \cite[Proposition 2.2]{E20} it was proved that a simple $3$-polytope $P\ne \Delta^3$ with $m$ faces is flag if and only if 
${\rm rk }\,H^4(\mathcal{Z}_P)=\frac{(m-2)(m-4)(m-6)}{3}$.

For simple polytopes of dimension $n\geqslant 4$ this is an open question (see \cite{FMW20}).
\begin{problem}
Is the property to be a flag $n$-polytope $B$-rigid in the class of simple $n$-polytopes?
\end{problem}
Also there is another known problem.
\begin{problem}
Is dimension $n$ of a simple polytope $B$-rigid in the class of all simple polytopes?
\end{problem}

It is easy to see that a polytope has no $4$-belts if and only if the multiplication
$$
H^3(\mathcal{Z}_P)\otimes H^3(\mathcal{Z}_P)\to H^6(\mathcal{Z}_P)
$$
is trivial.
Thus, the property to be a Pogorelov polytope is also $B$-rigid.

\begin{proposition} \cite{FMW15} (see also a proof in \cite[Proposition 8.24]{BE17S}) 
The subgroup $\boldsymbol{B}_k\subset H^{k+2}(\mathcal{Z}_P)$ is $B$-rigid for 
$4\leqslant k\leqslant m-2$. In particular, $\rk \boldsymbol{B}_k=\#\{\text{$k$-belts}\}$ for  $k\geqslant 4$ is a $B$-rigid
number.
\end{proposition}

In \cite{E20} it was proved that the property to be an almost Pogorelov polytope is $B$-rigid.
Namely, the polytopes
$I^3$ and $M_5\times I$ are $B$-rigid, since these are unique flag polytopes with $m=6$ and $m=7$ respectively.
Let $A_3(P)$ be a subgroup in $H^3(\mathcal{Z}_P)$, which is a kernel of the bilinear form defined
by the product map $H^3(\mathcal{Z}_P)\otimes H^3(\mathcal{Z}_P)\to H^6(\mathcal{Z}_P)$, that is
$$
A_3=\{x\in H^3(\mathcal{Z}_P)\colon x\cdot y=0\text{ for all } y\in H^3(\mathcal{Z}_P)\}.
$$
It is generated by elements $\widetilde{\{p,q\}}$ such that the faces $F_p$ and $F_q$ are not adjacent  
and do not belong to any $4$-belt $\mathcal{B}_4$.  

Denote by $I_7$ the image  of the mapping
$$
H^3(\mathcal{Z}_P)\otimes H^4(\mathcal{Z}_P)\to H^7(\mathcal{Z}_P)
$$
\begin{theorem}\cite[Theorem 3.12]{E20}
A flag simple $3$-polytope $P$ belongs to $\mathcal{P}_{aPog}\setminus\{I^3,M_5\times I\}$ if and only if 
$$
2\,\rk \boldsymbol{B}_4= \rk \left[H^3(\mathcal{Z}_P)/A_3(P)\right], \text{ and }\,\rk I_7=\rk \boldsymbol{B}_5+(m-5)\rk\,\boldsymbol{B}_4.
$$
In particular, the property to be an almost Pogorelov polytope is $B$-rigid. 
\end{theorem}
An almost Pogorelov polytope is ideal if and only if any its vertex lies on a unique quadrangle. This
is equivalent to the fact that $4p_4=f_0$, where $f_0$ is the number of vertices, and $p_4$ is the number of quadrangles.
For any simple $3$-polytope we have $f_0=2(m-2)$. For an almost Pogorelov polytope $P\ne I^3,M_5\times I$ 
we have $p_4=\rk\boldsymbol{B}_4$.
\begin{corollary}
A polytope in $\mathcal{P}_{aPog}\setminus\{I^3,M_5\times I\}$
belongs to $\mathcal{P}_{IPog}$ if and only if $2\,\rk\boldsymbol{B}_4=m-2$.
In particular, the property to be an ideal almost Pogorelov polytope is $B$-rigid.
\end{corollary}
\begin{problem}
We call by an {\it almost flag} polytope a simple $3$-polytope such that any its $3$-belt is trivial. Is the property to
be an almost flag polytope $B$-rigid? 
\end{problem}

\section{Three-dimensional classes}\label{3dsec}
In \cite{FMW15} it was proved that the set of elements 
$$
\{\pm \widetilde{\omega}\colon \omega\in N_2(P)\}\subset H^3(\mathcal{Z}_P)
$$ 
is $B$-rigid in the class of Pogorelov polytopes. 

The proof was based on the {\it annihilator lemma}. 
An {\em annihilator} of an element $r$ in a ring $R$ is 
${\rm Ann}_R(r)=\{s\in R\colon rs=0\}$.

\begin{lemma}[annihilator lemma, \cite{FMW15}]
Let $P\in \mathcal{P}_{Pog}$, and  $\alpha\in H=H^*(\mathcal{Z}_P,\mathbb Q)$: 
$$
\alpha=\sum\limits_{\omega\in N_2(P)}r_{\omega}\widetilde\omega\quad\text{with }|\{\omega\colon r_{\omega}\ne 0\}|\geqslant 2
$$ 
Then  
$\dim {\rm Ann}_H (\alpha)<\dim {\rm Ann}_H (\widetilde \omega),\text { if }r_{\omega}\ne 0$.
\end{lemma}
For almost Pogorelov polytopes the annihilator lemma 
is not valid. A counterexample arises already for $P=As^3$, see \cite[Proposition 8.3]{E20}. 

\begin{definition}\label{good-bad}
Let $P$ be a simple $3$-polytope. An element $\omega'=\{s,t\}\in N_2(P)$ is {\it good} for an element 
$\omega=\{p,q\}\in N_2(P)$, if there is an $l$-belt $\mathcal{B}_l$ containing $F_s$ and $F_t$ such that  
either $F_p$ or $F_q$ does not belong to $\mathcal{B}_l$ and does not intersect at least one of the two connected components of $\mathcal{B}_l\setminus\{F_s,F_t\}$. 
\end{definition}

\begin{lemma}\cite[Lemma 6.1]{E20}
Let  $P$ be a simple $3$-polytope and $\alpha\in H=H^*(\mathcal{Z}_P,\mathbb Q)$:
$$
\alpha=\sum\limits_{\omega\in N_2(P)}r_{\omega}\widetilde\omega\quad\text{with }|\{\omega\colon r_{\omega}\ne 0\}|\geqslant 2.
$$ 
Then for any $\omega=\{p,q\}$ with  $r_{\omega}\ne 0$ we have 
\begin{enumerate}
\item $\dim {\rm Ann}_H (\alpha)\leqslant \dim {\rm Ann}_H (\widetilde \omega)$; 
\item $\dim {\rm Ann}_H (\alpha)< \dim {\rm Ann}_H (\widetilde \omega)$, if there is $\omega'=\{s,t\}$ 
with $r_{\omega'}\ne\varnothing$, which is good for $\omega$.
\end{enumerate}
\end{lemma}
\begin{corollary}\cite[Corollary 6.4]{E20}
Let $P$ be a simple $3$-polytope. If any $\omega'\in N_2(P)\setminus\{\omega\}$ is good for an element $\omega\in N_2(P)$, 
then for any isomorphism of graded rings $\varphi\colon H^*(\mathcal{Z}_P)\to H^*(\mathcal{Z}_Q)$ for a simple $3$-polytope $Q$
we have 
$\varphi(\widetilde{\omega})=\pm \widetilde{\omega'}$ for some $\omega'\in N_2(Q)$.
\end{corollary}
It follows from the SCC for an almost Pogorelov polytope $P$ that for $\omega=\{p,q\}$, where $F_p$ and $F_q$ are quadrangles, 
any other element $\omega'\in N_2(P)$ is good. Thus, for any isomorphism of graded rings 
$\varphi\colon H^*(\mathcal{Z}_P)\to H^*(\mathcal{Z}_Q)$ for a simple $3$-polytope $Q$
we have 
$\varphi(\widetilde{\omega})=\pm \widetilde{\omega'}$ for some $\omega'\in N_2(Q)$. To prove that $\omega'$ also
corresponds to a pair of quadrangles an additional technique is needed.

First, for any polytope $P\in \mathcal{P}_{aPog}\setminus\{I^3,M_5\times I\}$ and any its $4$-belt
$\mathcal{B}_4$ define a subgroup $G(\mathcal{B}_4)$ in $\boldsymbol{H}_3=H^3(\mathcal{Z}_P)/A_3$ 
generated by two cosets $\widetilde{\{p,q\}}+A_3(P)$ corresponding to pairs of opposite faces of $\mathcal{B}_4$.

\begin{lemma} \cite[Lemma 7.1]{E20} 
The collection of subgroups 
$$
\{G(\mathcal{B}_4)\colon \mathcal{B}_4\text{ is a $4$-belt}\}
$$
is $B$-rigid in the class $\mathcal{P}_{aPog}\setminus\{I^3,M_5\times I\}$.
\end{lemma}

As a corollary the set  $\{\pm \widetilde{\mathcal{B}_4}\}$ of generators of $\boldsymbol{B}_4$ is $B$-rigid
in the class $\mathcal{P}_{aPog}\setminus\{I^3,M_5\times I\}$, since these are generators of the images of the mappings
$G(\mathcal{B}_4)\otimes G(\mathcal{B}_4)\to H^6(\mathcal{Z}_P)$. Therefore, for 
$P,Q\in \mathcal{P}_{aPog}\setminus\{I^3,M_5\times I\}$
any isomorphism of graded rings $H^*(\mathcal{Z}_P)\to H^*(\mathcal{Z}_Q)$ induces a bijection
$\varphi_4$ between the sets of quadrangles of  $P$ and $Q$ by the rule
$$
\varphi_4(F_i)=F_{i'}', \text{ where $\varphi(\widetilde{\mathcal{B}_i})=\pm \widetilde{\mathcal{B}_{i'}'}$ for 
$4$-belts  $\mathcal{B}_i$ and $\mathcal{B}_{i'}'$ around $F_i$ and $F_{i'}'$.}
$$

On the base of SCC for almost Pogorelov polytopes the next result was proved.
\begin{lemma} \cite[Lemma 7.5]{E20}
The collection of cosets of the form
$$
\pm\widetilde{\{p,q\}}+A_3(P), 
$$
where $F_p$ and $F_q$ are opposite faces of a $4$-belt, is $B$-rigid in the class 
$\mathcal{P}_{aPog}\setminus\{I^3,M_5\times I\}$.
\end{lemma}

The next important step is the following 
\begin{lemma}\cite[Proposition 7.9]{E20}\label{44lemma}
Let $P,Q\in \mathcal{P}_{aPog}\setminus\{I^3,M_5\times I\}$, and let $F_p$ be a quadrangle of $P$ not adjacent to a face $F_q$.
Assume that for an isomorphism of graded rings $\varphi\colon H^*(\mathcal{Z}_P)\to H^*(\mathcal{Z}_Q)$ we have  
$\varphi(\widetilde{\{p,q\}})=\pm \widetilde{\{s,t\}}$ for some $\{s,t\}\in N_2(Q)$. Then $p'\in \{s,t\}$ for
$F_{p'}'=\varphi_4(F_p)$. In particular, $\varphi(\widetilde{\{p,q\}})=\pm \widetilde{\{p',q'\}}$ for quadrangles $F_p$ and $F_q$, and the set 
$$
\{\pm\widetilde{\{p,q\}}\colon \text{ $F_p$ and $F_q$ are quadrangles}\}\subset H^3(\mathcal{Z}_P)
$$
is $B$-rigid in the class $\mathcal{P}_{aPog}\setminus\{I^3,M_5\times I\}$.
\end{lemma}
The proof is based on the following combinatorial fact and lemmas.
\begin{lemma}\cite[Lemma 7.11]{E20}
Let $P\in \mathcal{P}_{aPog}\setminus\{I^3,M_5\times I\}$. 
Then for any three pairwise different faces $\{F_i,F_j,F_k\}$ such that $F_i\cap F_j=\varnothing$, 
$F_k$ is a quadrangle, and at least one of the faces $F_i$ and $F_j$ is not adjacent to
$F_k$, there exists an $l$-belt $\mathcal{B}_l$ such that $F_i,F_j\in \mathcal{B}_l$, 
$F_k\notin \mathcal{B}_l$,  and $\mathcal{B}_l$ does not contain any of the two pairs of opposite faces of 
the $4$-belt around $F_k$.
\end{lemma}
\begin{lemma}\cite[Lemma 7.8]{E20}
Let $P\in \mathcal{P}_{aPog}\setminus\{I^3,M_5\times I\}$, and let $\mathcal{B}_k$ be a $k$-belt passing through a quadrangle 
$F_i$ and its adjacent faces $F_p$ an $F_q$. Let $x$ be a generator of $\widetilde{H}^1(P_{\tau})=\mathbb Z$
for $\tau=\omega(\mathcal{B}_k)$ (in this case $x=\pm\widetilde{\mathcal{B}_k}$), or $\tau=\omega(\mathcal{B}_k)\sqcup \{r\}$,
where $F_r$ either is not adjacent to faces in $\mathcal{B}_k$, or is adjacent to exactly one face in $\mathcal{B}_k$.
Then $x$ is divisible by any element in the coset $\widetilde{\{p,q\}}+A_3(P)$.
\end{lemma}
\begin{lemma}\cite[Corollary 7.7]{E20}\label{pqlemma}
Let $P,Q\in \mathcal{P}_{aPog}\setminus\{I^3,M_5\times I\}$, and let $\mathcal{B}_k$ be a $k$-belt 
passing through a quadrangle $F_i$ and its adjacent faces $F_p$ and $F_q$. 
Then for any isomorphism of graded rings $\varphi\colon H^*(\mathcal{Z}_P)\to H^*(\mathcal{Z}_Q)$ we have
$$
\varphi(\widetilde{\mathcal{B}_k})=\sum_j\mu_j\widetilde{\mathcal{B}_{k,j}'}
$$ 
for $k$-belts $\mathcal{B}_{k,j}'$ of $Q$ such that  for any $\mu_j\ne 0$ the belt $\widetilde{\mathcal{B}_{k,j}'}$ 
passes through the non-adjacent faces $F_{p'}'$ and $F_{q'}'$ of $Q$,
where $\varphi(\widetilde{\{p,q\}}+A_3(P))=\pm \widetilde{\{p',q'\}}+A_3(Q)$. Moreover, $F_{p'}'$ and $F_{q'}'$ are adjacent to 
the quadrangle $F_{i'}'=\varphi_4(F_i)$ of $Q$.
\end{lemma}

\section{Classes corresponding to belts}\label{Bcsec}
On the base of the $B$-rigid set of elements corresponding to pairs of non-adjacent faces
in \cite{FMW15} it was proved that the set of generators of $\boldsymbol{B}_k$:
$$
\{\pm \widetilde{\mathcal{B}_k}\colon \mathcal{B}_k - \text{ a $k$-belt}\}\subset H^{k+2}(\mathcal{Z}_P)
$$ 
is $B$-rigid in the class of Pogorelov polytopes. 

As we mentioned above for the class $\mathcal{P}_{aPog}\setminus\{I^3,M_5\times I\}$ this result
is valid for $k=4$. Also it can be proved \cite[Corollary 6.8]{E20} that if 
$P\in\mathcal{P}_{aPog}\setminus\{I^3,M_5\times I\}$ and a $k$-belt $\mathcal{B}_k$ does not have 
common points with quadrangles, then under any isomorphism of graded rings $H^*(\mathcal{Z}_P)\to H^*(\mathcal{Z}_Q)$ 
for a simple $3$-polytope $Q$ the element $\widetilde{\mathcal{B}_k}$ is mapped to $\pm\widetilde{\mathcal{B}_k'}$ for 
some $k$-belt $\mathcal{B}_k'$ of $Q$. 

For an ideal almost Pogorelov polytope any face, which is not a quadrangle, is surrounded by a $(2k)$-belt containing $k$ quadrangles, $k\geqslant3$. 

\begin{lemma}\cite[Lemma 7.12]{E20}
Let $P,Q\in \mathcal{P}_{aPog}\setminus\{I^3,M_5\times I\}$, and let $\mathcal{B}_{2k}$ be a $(2k)$-belt of $P$ containing 
$k$ quadrangles. Then for any isomorphism of graded rings $\varphi\colon H^*(\mathcal{Z}_P)\to H^*(\mathcal{Z}_Q)$  
we have  $\varphi(\widetilde{\mathcal{B}_{2k}})=\pm \widetilde{\mathcal{B}_{2k}'}$ for a $(2k)$-belt $\mathcal{B}_{2k}'$ in $Q$
containing $k$ quadrangles. In particular, the set of elements 
$$
\{\pm \widetilde{\mathcal{B}_{2k}}\colon \mathcal{B}_{2k} - \text{a $(2k)$-belt containing $k$ quadrangles}\}\subset H^{2k+2}(\mathcal{Z}_P)
$$ 
is $B$-rigid in the class $\mathcal{P}_{aPog}\setminus\{I^3,M_5\times I\}$.
\end{lemma}
The proof is based on Lemmas \ref{44lemma} and \ref{pqlemma}, and 
\begin{lemma}\cite[Corollary 6.10]{E20}
Let $P$ be a simple $3$-polytope and $\omega=\{p,q\}\in N_2(P)$. Then an element
$x=\sum_j\mu_j\widetilde{\mathcal{B}_{k,j}}\in\boldsymbol{B}_k$ is divisible by $\widetilde{\omega}$
if and only if each $\widetilde{\mathcal{B}_{k,j}}$ with $\mu_j\ne 0$ is divisible by $\widetilde{\omega}$,
and if and only if each belt $\mathcal{B}_{k,j}$ with $\mu_j\ne 0$ contains $F_p$ and $F_q$.
\end{lemma}

\section{Classes corresponding to trivial belts}\label{Btsec}

Then in \cite{FMW20} it was proved that the set of elements 
$$
\{\pm \widetilde{\mathcal{B}_k}\colon \mathcal{B}_k - \text{ a $k$-belt around a face}\}\subset H^{k+2}(\mathcal{Z}_P)
$$ 
is $B$-rigid in the class of Pogorelov polytopes. 
This induces a bijection between the sets of faces of polytopes.

We have an analog of this property.

\begin{lemma}\cite[Lemma 7.13]{E20}\label{2ktrivlemma}
Let $P,Q\in \mathcal{P}_{aPog}\setminus\{I^3,M_5\times I\}$, and let $\mathcal{B}_{2k}$ be a trivial $(2k)$-belt of $P$ containing 
$k$ quadrangles. Then for any isomorphism of graded rings $\varphi\colon H^*(\mathcal{Z}_P)\to H^*(\mathcal{Z}_Q)$  
we have  $\varphi(\widetilde{\mathcal{B}_{2k}})=\pm \widetilde{\mathcal{B}_{2k}'}$ for a trivial $(2k)$-belt 
$\mathcal{B}_{2k}'$ in $Q$ containing $k$ quadrangles. In particular, the set of elements 
$$
\{\pm \widetilde{\mathcal{B}_{2k}}\colon \mathcal{B}_{2k} - \text{a trivial $(2k)$-belt containing $k$ quadrangles}\}\subset H^{2k+2}(\mathcal{Z}_P)
$$ 
is $B$-rigid in the class $\mathcal{P}_{aPog}\setminus\{I^3,M_5\times I\}$.
\end{lemma}
The idea of the proof is similar to the idea of the proof for Pogorelov polytopes from \cite{FMW15}, 
namely, we consider elements in $H^{k+3}(\mathcal{Z}_P)$ having common divisors of a proper form
with $\widetilde{\mathcal{B}_{2k}}$. For Pogorelov polytopes it is sufficient to consider elements $\widetilde{\omega}$,
$\omega\in N_2(P)$.
But for almost Pogorelov polytopes we need to consider classes corresponding to pairs of quadrangles, and  
cosets of the form $\pm\widetilde{\{p,q\}}+A_3(P)$ corresponding to faces $F_p$ and $F_q$ adjacent to a quadrangle in the belt.

For ideal almost Pogorelov polytopes Lemma \ref{2ktrivlemma}
induces a bijection between the sets of faces of $P$ and $Q$, which are not quadrangles.
Together with $\varphi_4$ this bijections form a bijection between the sets of faces of $P$ and $Q$, since
any trivial belt of a polytope  $P\in\mathcal{P}_{aPog}\setminus\{I^3,M_5\times I\}$ surrounds exactly one face (for otherwise, 
$P$ is a prism and has adjacent quadrangles).

\section{Belts around adjacent faces}\label{Basec}
Then in \cite{FMW15} it was proved that for Pogorelov polytopes the images of adjacent faces under the mapping induced
by the isomorphism of graded rings $\varphi\colon H^*(\mathcal{Z}_P)\to H^*(\mathcal{Z}_Q)$ are adjacent. 
This follows from the fact that the trivial belts $\mathcal{B}_1$ and  $\mathcal{B}_2$ surround adjacent faces if and only if 
$\widetilde{\mathcal{B}_1}$ and $\widetilde{\mathcal{B}_2}$ have exactly one common divisor among $\widetilde{\omega}$, 
$\omega\in N_2(P)$. This finishes the proof of $B$-rigidity of any Pogorelov polytope. 

For ideal almost Pogorelov polytopes we have analogous facts.
\begin{lemma}\cite[Lemma 7.14]{E20}
Let $P\in \mathcal{P}_{aPog}\setminus\{I^3,M_5\times I\}$, $F_i$ be its quadrangle, and let 
$\mathcal{B}_{2k}$ be the $(2k)$-belt containing $k$ quadrangles and surrounding a face $F_j$. 
Then $F_i$ and $F_j$ are adjacent if and only if the element $\widetilde{\mathcal{B}_{2k}}$ is divisible by any 
element in the coset $\widetilde{\{p,q\}}+A_3(P)$ for one of the two pairs of opposite  faces $\{F_p, F_q\}$ of the $4$-belt
$\mathcal{B}_i$ around $F_i$. In particular,  if $F_i$ and $F_j$ are adjacent, then for any isomorphism of graded rings 
$\varphi\colon H^*(\mathcal{Z}_P)\to H^*(\mathcal{Z}_Q)$  for a polytope $Q\in \mathcal{P}_{aPog}\setminus\{I^3,M_5\times I\}$ 
we have  $\varphi(\widetilde{\mathcal{B}_{2k}})=\pm \widetilde{\mathcal{B}_{2k}'}$ for a $(2k)$-belt 
$\mathcal{B}_{2k}'$ containing $k$ quadrangles and surrounding a face $F_{j'}'$ adjacent to 
$F_{i'}'=\varphi_4(F_i)$.
\end{lemma}

\begin{lemma}\cite[Lemma 7.15]{E20}
Let $P\in \mathcal{P}_{aPog}\setminus\{I^3,M_5\times I\}$, and let 
$\mathcal{B}_{2k}$ be the $(2k)$-belt containing $k$ quadrangles and surrounding a face $F_i$, and
$\mathcal{B}_{2l}$ be the $(2l)$-belt containing $l$ quadrangles and surrounding a face $F_j$. 
Then $F_i$ and $F_j$ are adjacent if and only if the elements $\widetilde{\mathcal{B}_{2k}}$  and
$\widetilde{\mathcal{B}_{2l}}$ have exactly one common divisor among elements
$\widetilde{\{p,q\}}$ corresponding to pairs of quadrangles. 
In particular,  if $F_i$ and $F_j$ are adjacent, then for any isomorphism of graded rings 
$\varphi\colon H^*(\mathcal{Z}_P)\to H^*(\mathcal{Z}_Q)$  for a polytope $Q\in \mathcal{P}_{aPog}\setminus\{I^3,M_5\times I\}$ 
we have  $\varphi(\widetilde{\mathcal{B}_{2k}})=\pm \widetilde{\mathcal{B}_{2k}'}$ for a $(2k)$-belt 
$\mathcal{B}_{2k}'$ containing $k$ quadrangles and surrounding a face $F_{i'}'$,  and 
$\varphi(\widetilde{\mathcal{B}_{2l}})=\pm \widetilde{\mathcal{B}_{2l}'}$ for a $(2l)$-belt 
$\mathcal{B}_{2l}'$ containing $l$ quadrangles and surrounding a face $F_{j'}'$ 
adjacent to $F_{i'}'$.
\end{lemma}

This finishes the proof of $B$-rigidity of any ideal almost Pogorelov polytope.

\section{New cohomologically rigid families of manifolds}\label{Newsec}
\begin{definition}
A family $\mathcal{M}$ of manifolds is called {\it cohomologically rigid} over the ring $R$ if two manifolds 
$M_1,M_2\in\mathcal{M}$ are diffeomorphic if and only if there is an isomorphism of graded rings 
$H^*(M_1,R)\simeq H^*(M_2,R)$.
\end{definition}
$B$-rigidity of any Pogorelov polytope gives rise to cohomologically rigid over $\mathbb Z$ family of moment-angle manifolds.  Moreover, all the theory in \cite{FMW15} works over any field $\mathbb F$, therefore this family is
cohomologically rigid over any field $\mathbb F$. 

If we choose for any Pogorelov polytope $P$ some quasitoric manifold $M^6(P)$, then the family $\{M^6(P)\}$ of $6$-dimensional
manifolds is cohomologically rigid over $\mathbb Z$ or any field $\mathbb F$. The problem is that there is no canonical
way to choose such a manifold. 
Also, if we choose for each Pogorelov polytope a small cover $R(P)$, then the family of $3$-dimensional manifolds
$\{R(P)\}$ is cohomologically rigid over $\mathbb Z_2$. 

Theorem \ref{C6th} implies that the families of $6$-dimensional quasitoric manifolds
and $3$-dimensional small covers corresponding to the same Pogorelov polytope $P$ are cohomologically rigid  over
$\mathbb Z$ and $\mathbb Z_2$ respectively. 

For ideal almost Pogorelov polytopes we also obtain a cohomologically rigid family of moment-angle manifolds. For 
a polytope $P$ with $m$ faces (recall that $m=2(p_4+1)$, where $p_4$ is the number of quadrangles) the manifold
$\mathcal{Z}_P$ has dimension $m+3=2p_4+5$. 

For ideal almost Pogorelov polytopes the advantage is that there is a canonical way to choose a characteristic function. 
Namely, it is known
that a $3$-polytope can be coloured in $3$ colours such that adjacent faces have different colours 
if and only if any its face has an even number of edges (see \cite{I01,J01}). Moreover,
such a colouring is unique up to a permutation of colours. This gives a canonical characteristic function 
$\boldsymbol{\Lambda}_P$ for $P$ up to a permutation of basis vectors in $\mathbb Z^3$ (this function appears already in \cite[Example 1.15(1)]{DJ91}). 
Since a change of basis in $\mathbb Z^3$ and $\mathbb Z_2^3$ does 
not change the (weak equivariant) diffeomorphic type of manifolds  
$M(P,\Lambda)$ and $R(P,\Lambda_2)$, for any ideal almost Pogorelov polytope $P$ there is a canonical
quasitoric manifold $\boldsymbol{M}^6(P)$ and a canonical small cover $\boldsymbol{R}^3(P)$. 
We obtain three families of manifolds.
\begin{enumerate}
\item The family $\mathcal{Z}_{IPog}$ of moment-angle manifolds corresponding to ideal almost Pogorelov polytopes. 
\item The family $\mathcal{M}_{IPog}$ of canonical $6$-dimensional quasitoric manifolds $\boldsymbol{M}^6(P)$ 
corresponding to ideal almost Pogorelov polytopes.
\item The family $\mathcal{R}_{IPog}$ of canonical $3$-dimensional small covers $\boldsymbol{R}^3(P)$ 
corresponding to ideal almost Pogorelov polytopes.
\end{enumerate}
In \cite{DJ91} the manifolds from the second and the third families are called {\it pullbacks from the linear model}.

Each family is parametrised by ideal almost Pogorelov polytopes. On the other hand, as mentioned in the introduction,
ideal almost Pogorelov polytopes are in bijection with pairs $\{Q,Q^*\}$, where $Q$ is a combinatorial combinatorial convex 
(not necessarily simple) $3$-polytope, and $Q$ is its dual polytope. Namely for a polytope $Q$ the ideal almost Pogorelov 
polytope $P$ is obtained by cutting off all the vertices and then cutting off all the "old"  edges. 
The triangulation  $\partial P^*$ is the barycentric subdivision of $\partial Q^*$. Since all the arguments in this paper 
and \cite{E20} work for any field $\mathbb F$ taken instead of the ring $\mathbb Z$, we obtain the following result.

\begin{theorem}\label{FCth}
\begin{enumerate}
\item The family $\mathcal{Z}_{IPog}$ of moment-angle manifolds $\{\mathcal{Z}_P\}$ 
corresponding to ideal almost Pogorelov polytopes is cohomologically rigid over $\mathbb Z$ or any field~$\mathbb F$.
\item  The family $\mathcal{M}_{IPog}$ of canonical $6$-dimensional quasitoric manifolds $\{\boldsymbol{M}^6(P)\}$ 
corresponding to ideal almost Pogorelov polytopes is cohomologically rigid over $\mathbb Z$ or any field~$\mathbb F$.
\item The family $\mathcal{R}_{IPog}$ of canonical $3$-dimensional small covers $\{\boldsymbol{R}^3(P)\}$ 
corresponding to ideal almost Pogorelov polytopes is cohomologically rigid over~$\mathbb Z_2$.
\end{enumerate}
\end{theorem}

\section{Hyperbolic geometry of $3$-dimensional manifolds corresponding to ideal almost Pogorelov polytopes}
\label{Hsec}
For introduction to hyperbolic geometry and Coxeter groups we recommend  \cite{VS88}.

First let us describe explicitly the small cover $\boldsymbol{R}^3(P)$ associated to an ideal almost Pogorelov polytope $P$.
Remind that 
$$
\boldsymbol{R}^3(P)=\mathbb Z_2^3\times P/\sim, 
$$
where $(t_1,p)\sim (t_2,q)$ if and only if  $p=q$,  and $t_1-t_2\in\langle\boldsymbol{\Lambda}_i\rangle_{F_i\ni p}$.
Here $\boldsymbol{\Lambda}_i=\boldsymbol{e}_{c(i)}$, where $c(i)$ is the colour of the face $F_i$. 
For convenience we will assume that quadrangles are coloured in the third colour.

Thus, $\boldsymbol{R}^3(P)$ is clued from $8$ copies of $P$ of the form $t\times P$, $t\in\mathbb Z_2^3$. 
Each point in $t\times P$ lying in relative interior of
\begin{itemize}
\item the polytope $P$ belongs only to one polytope $t\times P$; 
\item the face $F_i$ belongs to two polytopes: $t\times P$ and $(t+\boldsymbol{\Lambda}_i)\times P$;
\item the edge $F_i\cap F_j$ belongs to four polytopes $t\times P$, $(t+\boldsymbol{\Lambda}_i)\times P$, 
$(t+\boldsymbol{\Lambda}_j)\times P$, and $(t+\boldsymbol{\Lambda}_i+\boldsymbol{\Lambda}_j)\times P$;
\item a vertex belongs to all the eight polytopes.
\end{itemize}

Let $Q$ be an ideal right-angled polytope such that $P$ is combinatorially obtained from $Q$ by cutting off vertices at infinity. 
We can realize $P$ as a part of $Q$ obtained by cutting off vertices  at infinity by small horospheres centered at them. The intersection of a horosphere with the polytope is a Euclidean rectangle.

In the polytope $P$ the faces intersecting any quadrangle by opposite edges have the same colour. 
When we pass one of these faces $F_i$, we move from $t\times P$ to $(t+\boldsymbol{\Lambda}_i)\times P$.
When we pass the other face $F_j$ in the new polytope  we move to 
$(t+\boldsymbol{\Lambda}_i+\boldsymbol{\Lambda}_j)\times P=t\times P$. Therefore,
in $\boldsymbol{R}(P)$ each quadrangle corresponds to a flat torus $T^2$ glued from $4$ copies
of a Eucledean rectangle). Since all the dihedral angles of $Q$ are right, outside this disjoint set of tori, the manifold 
has a hyperbolic structure. When we shrink the horospheres to points, then the tori  are also shrinked to points, 
and the manifold $\boldsymbol{R}(P)$ is transformed into a space glued from two parts $\boldsymbol{R}_1$ and 
$\boldsymbol{R}_2$ along finite sets of 
points corresponding to ideal vertices of $Q$.  If we delete these points, 
then each part $\widehat{\boldsymbol{R}_i}$ is a hyperbolic manifold of finite volume and is glued of four copies of $Q$:
$t\times Q$, $(t+\boldsymbol{e}_1)\times Q$, 
$(t+\boldsymbol{e}_2)\times Q$, and $(t+\boldsymbol{e}_1+\boldsymbol{e}_2)\times Q$. It
corresponds to a part of $\boldsymbol{R}(P)$ bounded by tori. Passing the quadrangles we 
add the vector $\boldsymbol{e}_3$ and move to the other part. Moreover, each part is a manifold
with boundary homotopy equivalent to $\widehat{\boldsymbol{R}_i}$, and $\boldsymbol{R}(P)$
is a double of this manifold: it is glued from two equal manifolds along boundaries.

The hyperbolic manifold $\widehat{\boldsymbol{R}_i}$ of finite volume can be obtained from $Q$
by a construction similar to Construction \ref{Vconstr} (see \cite[Section 5]{V17}). 
Namely, let $Q$ be an ideal right-angled polytope in the Lobachevsky space $\mathbb L^3$. 
There is a right-angled Coxeter  group $G(Q)$ generated by reflections in faces of $Q$. It it known that 
$$
G(Q)=\langle \rho_1,\dots,\rho_m\rangle/(\rho_1^2,\dots,\rho_m^2, \rho_i\rho_j=\rho_j\rho_i
\text{ for all }F_i\cap F_j\ne\varnothing),
$$
where $\rho_i$ is the reflection in the face $F_i$.

The colouring of $P$ into three colours induces the colouring of $Q$ into $2$ colours (black and white).

The mapping $\boldsymbol{\Lambda}$ defines a homomorphism $\varphi\colon G(Q)\to \mathbb Z_2^2=
\langle \boldsymbol{e}_1,\boldsymbol{e}_2\rangle\subset\mathbb Z_2^3=
\langle \boldsymbol{e}_1,\boldsymbol{e}_2,\boldsymbol{e}_2\rangle$ by the 
rule $\varphi(\rho_i)=\boldsymbol{\Lambda}_i$. It can be shown that $\Ker\varphi$ acts on $\mathbb L^3$ freely
(this follows from the fact that the stabiliser of a point in $Q$ by the action of $G(Q)$
is generated by reflections in faces containing this point). Then $\mathbb L^3/\Ker\varphi$ is a hyperbolic
manifold of finite volume homeomorphic to
$$
\mathbb Z_2^2\times Q/\sim, \text{ where }(t_1,p)\sim (t_2,q)\text{ if and only if } p=q, \text{ and } t_1-t_2\in\langle\boldsymbol{\Lambda}_i\rangle_{F_i\ni p}.
$$ 
A homeomorphism is given by the mapping $(t,p)\to \varphi^{-1}(t)\cdot p$. We see that 
this manifold is homeomorphic to $\widehat{\boldsymbol{R}_i}$.   
Moreover, the manifolds $\widehat{\boldsymbol{R}_1}$ and $\widehat{\boldsymbol{R}_2}$ correspond 
to two embeddings $\mathbb Z_2^2\to \mathbb Z_2^3$: $(x_1,x_2)\to (x_1,x_2,0)$ and $(x_1,x_2)\to (x_1,x_2,1)$.

\section{Acknowledgements}
The  author is grateful to Victor Buchstaber for his encouraging support and attention to this work, 
and Taras Panov and Alexander Gaifullin for useful discussions.

\end{document}